\documentclass[12pt]{amsart}

\title{Relative Inverse Limit Perfection of Derived Commutative Rings}
\author{Daniel Fink}
\address{Institut f\"ur Mathematik, Johannes Gutenberg-Universit\"at Mainz, 55099 Mainz, Germany}
\email{daniel.fink@uni-mainz.de}
\usepackage[pdftitle={lperfectfactorization},
pdfsubject={article},
pdfauthor={Daniel Fink},
pdfborder={0 0 0}]{hyperref}

\usepackage[margin=2.5cm]{geometry}
\usepackage{times}
\usepackage[utf8]{inputenc}
\usepackage[T1]{fontenc}
\usepackage{amsmath, amssymb, amsthm, mathtools}

\usepackage{tikz}
\usetikzlibrary{cd,babel,graphs}

\theoremstyle{definition}
\newtheorem{defi}{Definition}[section]
\newtheorem{ex}[defi]{Example}
\newtheorem{con}[defi]{Construction}
\newtheorem{rem}[defi]{Remark}

\theoremstyle{plain}
\newtheorem{thm}[defi]{Theorem}
\newtheorem{prop}[defi]{Proposition}
\newtheorem{lem}[defi]{Lemma}
\newtheorem{cor}[defi]{Corollary}

\newcommand{\NN}{\mathbf{N}}  
\newcommand{\ZZ}{\mathbf{Z}}

\newcommand{\FF}{\mathbf{F}}
\newcommand{\sslash}{\mathbin{/\mkern-6mu /}}

\newcommand{\Hom}{\operatorname{Hom}}
\newcommand{\fun}{\operatorname{Fun}}

\newcommand{\colim}{\operatorname{colim}}

\newcommand{\id}{\operatorname{id}}

\newcommand{\perfection}{\mathrm{perf}}

\newcommand{\Mod}{\operatorname{Mod}}
\newcommand{\dalg}{\operatorname{DAlg}}
\newcommand{\dalgperf}[1]{\operatorname{DAlg}^{\operatorname{perf}}_{\substack{#1}}}
\newcommand{\ACR}{\operatorname{ACR}}
\newcommand{\perf}[1]{\operatorname{ACR}^{\operatorname{perf}}_{\substack{#1}}}
\newcommand{\sperf}[1]{\operatorname{ACR}^{\operatorname{semi}}_{\substack{#1}}}

\renewcommand{\to}[1][]{\xrightarrow{\ #1\ }}
\newcommand{\ot}[1][]{\xleftarrow{\ #1\ }}

\newcommand{\into}[1][]{\xhookrightarrow{\ #1\ }}

\usepackage[backend=biber,style=numeric, natbib=true,hyperref=true,language=auto, autolang=other,maxnames=50, sortcase=false, minalphanames=4, maxalphanames=4, sorting=nyt]{biblatex}

\begin{document}

\begin{abstract}
We study the relative Frobenius map associated with a map of derived commutative rings over a field of positive characteristic. As part of this, we examine a relative analog of perfectness and construct a relative inverse limit perfection which, under suitable conditions on the base, serves as a right adjoint to the inclusion of relatively perfect algebras into the category of all algebras. 
Specializing to animated rings, we investigate relative versions of semiperfectness and \textit{F}-finiteness, and use these to show that any map of \textit{F}-finite animated rings factors into a free map of finite type, followed by a relatively perfect map, followed by a surjective map. We also show that, for a morphism of Noetherian \textit{F}-finite rings, the vanishing of the cotangent complex implies that the morphism is relatively perfect.
\end{abstract}

\maketitle


\section{Introduction}
\renewcommand{\thedefi}{\Alph{defi}}
A homomorphism of finite type between commutative rings is given by finitely many generators, subject to (possibly infinitely many) relations. The main goal of this article is to show that --- up to a relatively perfect map --- a similar result holds for maps of commutative \textit{F}-finite rings over the prime field $\FF_p$. 

To outline the construction, we consider relative analogs of \textit{F}-finiteness, semiperfectness, and perfectness, which specialize to their absolute counterparts when the base is $\FF_p$. Let $R\to S$ be a morphism of $\FF_p$-algebras, and consider the associated relative Frobenius homomorphism:
\[
S\otimes_{R,F} R \to S, \quad s\otimes r\mapsto s^pr.
\]
We call the $R$-algebra $S$ \emph{relatively \textit{F}-finite} if the relative Frobenius morphism is finite, \emph{relatively semiperfect} if it is surjective, and \emph{relatively perfect} if it is an isomorphism and $S$ and $F_*R$ are Tor-independent over $R$. 

Explicitly, relative \textit{F}-finiteness means that $S$ is finitely generated over its subring $S^p[R]$, while relative semiperfectness simply means $S^p[R] = S$. In the relative \textit{F}-finite case, choosing generators $s_1,\dots,s_n$ of $S$ over $S^p[R]$ yields a map from the polynomial ring $R[x_1,\dots,x_n]$ to $S$.
Over this new base, $S$ is relatively semiperfect, as $S^p[R,s_1,\dots,s_n]=S$. This gives the first part of the factorization.

Now, given a map $R\to S$ that makes $S$ relatively semiperfect over $R$, we aim to construct a relatively perfect $R$-algebra $T$ together with a surjective map $T\to S$ of $R$-algebras. If $R$ is perfect, this is done by computing the inverse limit perfection of $S$. For more general base rings, however, this perfection is too coarse and typically does not yield an $R$-algebra. Instead, we iterate the relative Frobenius map to form an inverse system:
\[
\dots\to[F\otimes \id_R] S\otimes_{R,F^3}R\to[F\otimes \id_R] S\otimes_{R,F^2}R\to[F\otimes \id_R] S\otimes_{R,F}R \to[F\otimes \id_R] S\otimes_{R}R\simeq S.
\]
Under suitable assumptions on $R$ --- for instance, if $R$ is regular Noetherian and \textit{F}-finite --- the inverse limit of this system gives exactly the desired relatively perfect algebra. The general case can then be reduced to this one by base change. When $R$ is a polynomial ring in finitely many variables over $\FF_p$, this agrees with the construction described in Remark~13.6 of Gabber’s article~\citep{Gab04}; see Corollary \ref{cor:introc} and Example \ref{ex:gabbers_construction_der} in the main text.

The preceding constructions extend naturally to the derived setting. We call a map of derived commutative $\FF_p$-algebras \emph{relatively perfect} if its derived relative Frobenius map is an equivalence. For maps of connective derived rings (animated rings), we define the relative notions of \textit{F}-finiteness and semiperfectness as before, using their $H^0$.

For details on \emph{derived commutative rings}, see the articles by Raksit~\citep{Rak20}\footnote{Raksit attributes the development of the theory to Mathew and also mentions an unpublished work-in-progress by Bhatt and Mathew; see \citep[\S 1.3]{Rak20}.} and Holeman~\citep{Hol23}.
For animated rings, see Lurie’s thesis~\citep{Lur04} and Spectral Algebraic Geometry \citep{Lur18}, where they are referred to as \emph{simplicial commutative rings}. See also the article of \v{C}esnavi\v{c}ius and Scholze~\citep{CS24}, where the term \emph{animated rings}, suggested by Clausen, was introduced.

\subsection*{Outline and Main Results}
In Sections \ref{sec:derivd_rel_frob} and \ref{sec:relative_perfectness}, we study the relative Frobenius map and the notion of relative perfectness for maps of general derived $\FF_p$-algebras. Using the usual notion of flatness for the underlying maps of $E_{\infty}$-rings, we obtain the following coconnectivity result (Lemma \ref{lem:flat_frobenius_implies_coconnective}).

\begin{lem}
\label{lem:introduction}
A derived $\FF_p$-algebra with a flat Frobenius endomorphism has no negative cohomology.
\end{lem}

In particular, animated $\FF_p$-algebras with a flat Frobenius endomorphism are discrete. Since flatness is preserved under base change, the property of having a flat absolute Frobenius map extends to relatively perfect algebras. Using this, we show that discrete relatively perfect algebras over a ring with a flat Frobenius endomorphism are Tor-independent (Proposition \ref{prop:rel_perf_tor-ind}). This, in turn, gives a base change result (Corollary~\ref{cor:perfect_base-change_der}) following Lemma~3.18 in the article of Bhatt and Scholze \citep{BS17}.

We also recall that the cotangent complex of a map of animated $\FF_p$-algebras is closely related to that of the associated relative Frobenius map (Lemma \ref{lem:cotangent_relative_frobenius}), which implies that the cotangent complex of a relatively perfect map vanishes (Proposition \ref{prop:cotangent_rel_perf_animated}). This provides a criterion for deciding whether a map of animated $\FF_p$-algebras is an equivalence (Corollary \ref{cor:equivalence_animated}).

We then associate to a map $R\to S$ of derived rings its \emph{relative Frobenius tower} (Construction \ref{con:derived_rel_frobenius_tower}), and show that if $(-)\otimes_{R,F}R$ commutes with limits, then the limit $S^{\perfection/R}$ of this tower is relatively perfect over $R$ (Lemma \ref{lem:derived_relative_perfection}) and serves as the universal relatively perfect $R$-algebra admitting a map to $S$ (Proposition \ref{prop:derived_relative_inverse_right_adjoint}). We call this the \emph{relative inverse limit perfection}.

In the remaining sections we consider maps $R\to S$ of animated commutative $\FF_p$-algebras. Section~\ref{sec:rel_semiperf} begins with a study of relative semiperfectness, as outlined above. This condition ensures that the relative inverse limit perfection $S^{\perfection/R}$ remains an animated ring and that the natural map $S^{\perfection/R}\to S$ is surjective on $H^0$ (Lemma \ref{lem:connective_surjective_der}). We first analyze several properties of the inverse limit perfection in the relative semiperfect case, and then turn to the study of \textit{F}-finite maps.

The main result will be established in Section \ref{sec:factorization}; see Theorem \ref{thm:F-finite_factorization_der} for the precise statement.

\begin{thm}
\label{thm:2intro2}
Any map $R \to S$ of \textit{F}-finite animated $\FF_p$-algebras admits a factorization
\[
R \to R^{\prime} \to T \to S
\]
where the first map is of finite type and free, the second is relatively perfect, and the third is surjective on $H^0$. Moreover, the explicit construction ensures:
\begin{enumerate}
\item[(i)] If $H^0(R)$ is Noetherian, and $R$ is $n$-coconnective, then $T$ is also $n$-coconnective.
\item[(ii)] If $R$ and $S$ are Noetherian, then $T$ is Noetherian as well.
\end{enumerate}
\end{thm}

The theorem has two interesting special cases. The first is when the base ring is $\mathbb{F}_p$.

\begin{cor}
\label{cor:introc}
For any \textit{F}-finite animated $\FF_p$-algebra $S$ there exists a free $\FF_p$-algebra of finite type $R^{\prime}$, and a relatively perfect $R^{\prime}$-algebra $T$ along with a map $T\to S$ that is surjective on $H^0$.
\end{cor}

The construction of $T$ follows the outline given earlier applied to the map $R^{\prime}\to H^0(S)$. For details, see Remark \ref{rem:f-compact}, and Lemmas \ref{lem:perfection_sees_only_H0_flat}, \ref{lem:derived_relative_perfection},   \ref{lem:connective_surjective_der}, and \ref{lem:f-finite_smooth_f-surjective}.

Since $T$ is relatively perfect over a polynomial ring in finitely many variables over $\FF_p$, it has remarkable properties: it is discrete, its Frobenius endomorphism makes its target a finite free module over its source (Remark \ref{rem:der_rel_perf_bc_prop}), and its cotangent complex over $\FF_p$ is equivalent to a finite free module concentrated in degree $0$ (Proposition  \ref{prop:cotangent_rel_perf_animated}). Moreover, if $H^0(S)$ is Noetherian, then $T$ is also Noetherian and thus regular (Proposition \ref{prop:perfection_is_Noetherian_der}).

The second special case is when $R$ and $S$ are discrete (Corollary \ref{cor:factorization_discrete_der} and Remark \ref{rem:final_remark}).

\begin{cor}
\label{thm:introduction_3}
Any map $R\to S$ of \textit{F}-finite $\FF_p$-algebras admits a factorization
\[
R \to R^{\prime} \to T \to S
\]
where $R^{\prime}$ is a polynomial ring over $R$ in finitely many variables, the (discrete) relative Frobenius map $F_{T/R^{\prime}}$ is an isomorphism, and the last map is surjective. 
If $R$ is Noetherian, then $R^{\prime} \to T$ is relatively perfect; if both $R$ and $S$ are Noetherian, then $T$ is Noetherian as well.
\end{cor}

Since the cotangent complex of a relatively perfect map vanishes, the factorization $R \to T \to S$ from the previous corollary, for maps of Noetherian \textit{F}-finite rings, is an $L$-smooth-by-surjective factorization in the sense of \citep{BF25}; the existence of this factorization was suggested by Bhatt.

Using the factorization results of Section \ref{sec:factorization}, we also prove that a map $R\to S$ of Noetherian $F$-finite $\FF_p$-algebras is relatively perfect if its cotangent complex $L_{S/R}$ vanishes (Theorem \ref{thm:acyclic_cotanget_implies_rel_perf}). Since $L_{S/R}\simeq 0$ is equivalent to formal \'etaleness for maps of Noetherian rings, we obtain the following result.

\begin{thm}
\label{thm:introduction4}
A map $R\to S$ of Noetherian $F$-finite $\FF_p$-algebras is relatively perfect if and only if it is formally \'etale. 
\end{thm}

From this we deduce that any formally smooth map $X\to Y$ of locally Noetherian \textit{F}-finite $\FF_p$-schemes admits, Zariski locally, a factorization $X\to \mathbb{A}^n_Y \to Y$ where the first morphism is formally \'etale and the second is the canonical projection (Corollary \ref{cor:factorizationofschemes}).

\subsection*{Notation and Terminology} 
We use the language of $\infty$-categories as developed by Lurie in \citep{Lur09,Lur17}.
We denote by $D(\FF_p)$ the derived $\infty$-category of $\FF_p$-vector spaces, equipped with its usual symmetric monoidal structure and $t$-structure. We regard $D(\mathbb{F}_p)$ as a derived algebraic context by taking $D(\mathbb{F}_p)^0$ to be the full subcategory of finite-dimensional $\mathbb{F}_p$-vector spaces concentrated in degree zero, following \cite[Definition 4.2.1]{Rak20}. By \citep[Theorem 7.1.2.13]{Lur17}, there is an equivalence $D(\FF_p)\simeq\Mod_{\FF_p}$, identifying this algebraic context with \citep[Example 4.3.1]{Rak20}, replacing $\ZZ$ by $\FF_p$.

We write 
$\operatorname{DAlg}_{\FF_p}=\operatorname{DAlg}(D(\FF_p))$ for the $\infty$-category of derived commutative\footnote{We always work in the commutative setting but often omit the adjective for brevity.} $\FF_p$-algebras \citep[Definition 4.2.22]{Rak20}. 

For $R\in \operatorname{DAlg}_{\FF_p}$ and $n\in \ZZ$, we denote by $H^n(R)\cong \Hom_{hD(\FF_p)}(\FF_p[-n],R)$ the $n$-th cohomology of the underlying object in $D(\FF_p)$. Here, $hD(\FF_p)$ denotes the homotopy category of $D(\FF_p)$, which is the classical derived category of $\FF_p$-vector spaces. By \citep[Proposition 4.2.27]{Rak20}, every derived $\FF_p$-algebra has an underlying $E_{\infty}$-algebra object in $D(\FF_p)$. In particular, for all $m,n\in \ZZ$, the multiplication map $R\otimes^L_{\FF_p} R \to R$ determines bilinear maps 
\[
H^m(R)\times H^n(R) \to \Hom_{hD(\FF_p)}\left(\FF_p[-m]\otimes^L_{\FF_p}\FF_p[-n], R\otimes^L_{\FF_p} R\right)\to H^{m+n}(R).
\]
These maps equip $H^{\star}(R)= \oplus_{n\in\ZZ}H^n(R)$ with the structure of a graded-commutative $\FF_p$-algebra, and $H^0(R)$ with the structure of an ordinary commutative $\FF_p$-algebra. We say that $R$ is \emph{discrete} if $R\simeq H^0(R)$, \emph{(cohomologically) $n$-connective} if $H^{k}(R)=0$ for all $k > n$; and \emph{(cohomologically) $n$-coconnective} if $H^{k}(R)=0$ for all $k< n$, where $n\in\ZZ$.
By \citep[Remark 4.3.3]{Rak20}, the $\infty$-category $\dalg_{\FF_p}\times_{D(\FF_p)}D(\FF_p)^{\le 0}$ of connective derived rings is canonically equivalent to the $\infty$-category $\ACR_{\FF_p}$ of animated (commutative) $\FF_p$-algebras.

If $R$ is a derived $\FF_p$-algebra, a derived $R$-algebra is a derived $\FF_p$-algebra $S$ together with a map $R\to S$ of derived $\FF_p$-algebras. We denote by $\dalg_R=\dalg_{\FF_p,R/}$ the $\infty$-category of derived $R$-algebras. If $R$ is connective, we denote by $\ACR_R=\ACR_{\FF_p,R/}$ the $\infty$-category of animated $R$-algebras.

For an ideal $I$ of a discrete ring $R$, we write $I^{[p^k]}$ for $F^k(I)R$, where $F=F_R$ denotes the Frobenius endomorphism of $R$. 

Given a map $R\to S$ of (derived) $\FF_p$-algebras, 
$S\otimes^L_{R,F}R$ denotes the pushout of $S\ot R \to[F]R$ in $\dalg_{\FF_p}$.

\subsection*{Acknowledgements}
The present article builds on and further develops ideas from \citep{BF25}.

I am grateful to \mbox{Manuel Blickle} for posing the question of a relative analog of the construction in \citep[Remark 13.6]{Gab04} and for encouraging me to investigate it. I also thank him for his continuous support, many helpful conversations, and for sharing an early manuscript of \citep{BBST}. I am also grateful to \mbox{Andreas Gieringer}, \mbox{Nikita M\"uller}, and \mbox{Luca Passolunghi} for insightful conversations on this topic, and especially to Luca for explaining to me the proofs of Lemma~\ref{lem:presentable} and Proposition~\ref{prop:derived_relative_inverse_right_adjoint}. 

The author acknowledges support from the Deutsche Forschungsgemeinschaft (DFG, German Research Foundation) through the Collaborative Research Centre TRR 326 \textit{Geometry and Arithmetic of Uniformized Structures}, project number 444845124.

\renewcommand{\thedefi}{\thesection.\arabic{defi}}

\section{The derived relative Frobenius map}
\label{sec:derivd_rel_frob}
We recall from Construction 2.4.1 in Holeman's article \citep{Hol23} the existence of a functor
\[
\dalg_{\FF_p} \to \fun(B\NN, \dalg_{\FF_p}),
\]
which assigns to each derived commutative $\FF_p$-algebra its Frobenius endomorphism, also referred to as the absolute Frobenius map. Here, $B\NN$ denotes the classifying simplicial set of the monoid $\NN$; see \cite[\href{https://kerodon.net/tag/04FY}{Subsection 04FY}]{kerodon}. As noted in \citep[Remark 2.4.2]{Hol23}, this functor preserves small limits and colimits, so the Frobenius endomorphism of any such limit or colimit is functorially determined by those of the rings in the diagram. Moreover, by \citep[Observation 2.2.19]{Hol23}, this construction specializes to the animated Frobenius when restricted to the subcategory of connective derived $\FF_p$-algebras.

The existence of the absolute Frobenius map enables us to define a relative version.

\begin{defi}
\label{def:relative_frobenius_der}
Let $R\to S$ be a map of derived $\FF_p$-algebras. The naturality square of the absolute Frobenius maps of $R$ and $S$ induce a commutative diagram:
\begin{equation}
\begin{tikzcd}
R \arrow[rr, "F"] \arrow[dd] &                                                               & R \arrow[dd] \arrow[ld] \\
                             & {S\otimes^L_{R,F}R} \arrow[rd, "F_{S/R}" description, dashed] &                         \\
S \arrow[rr, "F"] \arrow[ru] &                                                               & S                      
\end{tikzcd}
\tag{\(\star\)}
\label{eq:relative_frobenius_der}
\end{equation}
We call $F_{S/R}$ a \emph{relative Frobenius map of $S$ over $R$}.
\end{defi}

\begin{rem}
\label{rem:animated_relative_frobenius_der}
Let $R\to S$ be a map of animated commutative $\FF_p$-algebras. By \citep[Corollary 2.13]{Mao24} and its proof, this map can be expressed as a sifted colimit of maps $\FF_p[X] \to \FF_p[X,Y]$, where $X$ and $Y$ are finite sets, and the map restricts to the identity on $X$. Since the Frobenius maps of $R$ and $S$ are functorially induced by the Frobenius maps at each stage of this sifted diagram, the relative Frobenius map of $S$ over $R$ is similarly induced by the relative Frobenius maps of the sifted diagram. In other words, the relative Frobenius map of connective derived rings arises through the animation of the classical relative Frobenius map.

Moreover, 
since $H^0$ preserves colimits when restricted to the full subcategory of animated $\FF_p$-algebras, it follows that $H^0(F_{S/R})$ is the discrete relative Frobenius map of $H^0(S)$ over $H^0(R)$.
\end{rem}

We record the following basic properties:

\begin{rem} 
\label{rem:basic_properties_general_der}
\begin{enumerate}
\item[(i)] Two composable maps $R \to S \to T$ of derived $\FF_p$-algebras induce a commutative diagram of derived $R$-algebras:
\[
\begin{tikzcd}
{S\otimes^L_{R,F}R} \arrow[r, "F_{S/R}"] \arrow[d]                          & S \arrow[d]                              &   \\
{T\otimes^L_{R,F}R} \arrow[r] \arrow[rr, "F_{T/R}" description, bend right] & {T\otimes^L_{S,F}S} \arrow[r, "F_{T/S}"] & T
\end{tikzcd}
\]
\item[(ii)] Let $R'\ot R \to S$ be a diagram of derived $\FF_p$-algebras, and let $S'=R'\otimes^L_R S$. Then, the canonical diagram is cocartesian:
\[
\begin{tikzcd}
{S\otimes^L_{R,F}R} \arrow[r] \arrow[d, "F_{S/R}"'] & {S'\otimes^L_{R',F}R'} \arrow[d, "F_{S'/R'}"] \\
S \arrow[r]                                       & S'                                         
\end{tikzcd}
\]

\end{enumerate}
\end{rem}

\begin{rem}
\label{rem:limit-colimit-frobenius_der} 
Let $R$ be a derived $\FF_p$-algebra, and let $S:I\to \operatorname{DAlg}_{R}$ be a small diagram of $R$-algebras. Then, the relative Frobenius map of the limit is equivalent to the composition 
\[
(\lim_{i\in I}S_i) \otimes^L_{R,F} R\to \lim_{i\in I}(S_i \otimes^L_{R,F}R) \to[\lim_{i\in I}F_{S_i/R}] \lim_{i\in I}S_i.
\]
Similarly, the relative Frobenius map of the colimit is equivalent to 
\[
(\colim_{i\in I}S_i) \otimes^L_{R,F} R\to[\simeq] \colim_{i\in I}(S_i \otimes^L_{R,F}R) \to[\colim_{i\in I}F_{S_i/R}] \colim_{i\in I}S_i.
\]
\end{rem}

Using the following definition of flatness for derived rings, which is similar to \citep[Definition 7.2.2.10]{Lur17} for the underlying $E_{\infty}$-rings, we obtain a lemma that will be important for several arguments and results later in this article.

\begin{defi}
\label{def:flatness}
We call a map $R\to S$ of derived rings \emph{flat}, if and only if $S$ is a flat module over $R$, meaning that $H^0(S)$ is a flat over $H^0(R)$, and for all $k\in \ZZ$, the canonical map 
\[
H^k(R)\otimes_{H^0(R)}H^0(S)\to H^k(S),
\]
which sends $x\otimes y$ to $(H^k f)(x)\cdot y$, is an isomorphism.
\end{defi}

\begin{lem}
\label{lem:flat_frobenius_implies_coconnective}
Let $R$ be a derived $\FF_p$-algebra. If the absolute Frobenius map of $R$ is flat, then $R$ is coconnective. 
\end{lem}
\begin{proof}
The proof of \citep[Proposition 11.6]{BS17} by Bhatt and Scholze shows that the absolute Frobenius map of a derived $\FF_p$-algebra induces the zero map on negative cohomology groups. The claim then follows from the definition of flatness.
\end{proof}

In particular, for an animated commutative ring $R$, the absolute Frobenius map $F_R$ is flat if and only if $H^0(F_R)=F_{H^0(R)}$ is flat and $R$ is discrete.

If the base ring has a flat Frobenius map, the relative Frobenius map also induces the zero map on negative cohomology groups.

\begin{lem}
\label{lem:rel_frob_zero_on_cohomology}
Let $R\to S$ be a map of derived $\FF_p$-algebras, and suppose that the absolute Frobenius map $F_R$ is flat. Then $H^k(F_{S/R})=0$ for all $k\le -1$. 
\end{lem}
\begin{proof}
Since the absolute Frobenius map of $R$ is flat, applying $H^k$ to the relative Frobenius map
\[
S\otimes_{R,F}^L R 
\to[F_S\otimes \id_R] 
S\otimes_R^L R
\]
yields
\[
 H^k(S\otimes^L_{R,F} R) \cong H^k(S)\otimes_{H^0(R),F} H^0(R)
\to[H^k(F_S)\otimes \id_{H^0(R)}] 
H^k(S)\otimes_{H^0(R)} H^0(R) \cong H^k(S),
\]
which is zero for $k\le -1$. 
\end{proof}

\begin{rem}
\label{rem:canonical_fibersequences_connective}
For any map $R\to S$ of animated $\FF_p$-algebras, the diagram (\ref{eq:relative_frobenius_der}) canonically induces two fiber sequences of cotangent complexes:
\[
\begin{array}{c@{\quad}c@{\quad}c}
    L_{S/R} \otimes^L_{S,F} S \to L_{S/R} \to L_{F_{S/R}} & \quad \text{and} \quad & L_{F_R} \otimes^L_R S \to L_{F_S} \to L_{F_{S/R}.}
\end{array}
\]
\end{rem}

From the first of these exact triangles, it follows that the cotangent complex of the relative Frobenius map is determined by the cotangent complex of the underlying map of animated rings.

\begin{lem}
\label{lem:cotangent_relative_frobenius}
Let $R\to S$ be a map of animated commutative $\FF_p$-algebras. Then, there exists a canonical equivalence
\[
L_{F_{S/R}}\simeq L_{S/R} \oplus (L_{S/R}\otimes^L_{S,F} S[1])
\]
of $S$-modules.
\end{lem}
\begin{proof}
We show that the map $L_{S/R} \otimes^L_{S,F} S \to L_{S/R}$ is zero. The stated equivalence then follows formally, since $L_{F_{S/R}}$ is the cofiber of this map. 

We can express $R\to S$ as a sifted colimit of maps of the form $\FF_p[X]\to \FF_p[X,Y]$, where $X$ and $Y$ are finite sets; see the paragraph preceding Definition 2.33 in Mao’s article \citep{Mao24}. The Frobenus map of $R\to S$ is obtained by passing to the colimit of the Frobenius morphisms in this system. Consequently, the above map between cotangent complexes is induced by the maps arising from the commutative diagrams
\[
\begin{tikzcd}
{\FF_p[X]} \arrow[r, "F"] \arrow[d] & {\FF_p[X]} \arrow[d] \\
{\FF_p[X,Y]} \arrow[r, "F"]           & {\FF_p[X,Y]}          
\end{tikzcd}
\]
in each stage of the sifted diagram. Thus, it suffices to prove the claim for $\FF_p[X]\to \FF_p[X,Y]$. This follows from a straightforward calculation: any element of the form $d(r)\otimes b$ is mapped to $d(r^p)b^p=0$.
\end{proof}

\section{Relative Perfectness} 
\label{sec:relative_perfectness}
We recall that a derived $\FF_p$-algebra is called \emph{perfect} if and only if its Frobenius endomorphism is an equivalence. In this section, we discuss a relative analog of perfectness using the relative Frobenius map.

\begin{defi}
A map $R\to S$ of derived commutative $\FF_p$-algebras is \emph{relatively perfect} if and only if the relative Frobenius map $F_{S/R}$ is an equivalence. In this case, we say that \emph{$S$ is a relatively perfect $R$-algebra}.

We denote by $\dalgperf{R}$ the full subcategory of $\dalg_R$ consisting of relatively perfect $R$-algebras.
\end{defi}

\begin{ex}
Any equivalence of derived $\FF_p$-algebras is relatively perfect. 

A derived $\FF_p$-algebra is relatively perfect over $\FF_p$ if and only if it is perfect. More generally, let $R\to S$ be a map of derived $\FF_p$-algebras with $R$ perfect. Then $F_S$ factors as an equivalence followed by $F_{S/R}$, so $S$ is relatively perfect over $R$ if and only if it is perfect itself.
\end{ex}

\begin{rem}
\label{rem:equivalence_of_spaces_der}
Let $R$ be a derived $\FF_p$-algebra, and let $T$ and $S$ be derived $R$-algebras. If $T$ is relatively perfect, postcomposition with $F_{S/R}$ induces an equivalence of mapping anima
\[
\operatorname{Map}_{\ACR_R}(T,S\otimes_{R,F}R)\to \operatorname{Map}_{\ACR_R}(T,S),
\]
with an inverse given by a composition
\[
\operatorname{Map}_{\ACR_R}(T,S) \to[(-)\otimes^L_{R,F}R] \operatorname{Map}_{\ACR_R}(T\otimes^L_{R,F}R,S\otimes^L_{R,F}R) \to[\circ F_{T/R}^{-1}] \operatorname{Map}_{\ACR_R}(T ,S\otimes^L_{R,F}R).
\]
\end{rem}

\begin{rem}
\label{rem:tor_independence_relative_perfect_der}
For a map $R\to S$ of discrete rings, our definition of relative perfectness requires that the classical (non-derived) relative Frobenius map is an isomorphism, and that $S$ and $F_*R$ are Tor-independent over $R$. If both rings are Noetherian, Tor-independence holds automatically whenever the classical relative Frobenius map is an isomorphism. Indeed, in this case the map $R\to S$ is formally smooth (even formally \'etale) by \citep[Théorème 21.2.7]{EGAIV} and thus flat by \citep[Théorème
19.7.1]{EGAIV}. If the base ring is regular Noetherian, our definition also coincides with the classical notion.
\end{rem}

\begin{rem}
\label{rem:der_rel_perf_bc_prop}
Let $R$ be a derived $\FF_p$-algebra. By definition, the absolute Frobenius map of any relatively perfect derived $R$-algebra inherits all properties of $F_R$ that are stable under base change. For instance, if $R$ has a flat Frobenius endomorphism, then any relatively perfect $R$-algebra also has a flat Frobenius endomorphism and is therefore coconnective by Lemma \ref{lem:flat_frobenius_implies_coconnective}. 
\end{rem}

\begin{ex}
\label{ex:etale_vs_rel_perf_der}
The notion of relative perfectness for discrete $\FF_p$-algebras is closely related to various forms of \'etaleness. We have the following implications:
\begin{align*}
\text{\'etale} 
&\implies \text{ind-\'etale} 
\implies \text{weakly \'etale} 
\implies \text{relatively perfect}
\\
&\implies \text{acyclic cotangent complex} 
\implies \text{formally \'etale};
\end{align*}
see \cite[\href{https://stacks.math.columbia.edu/tag/097N}{Lemma 097N}]{stacks-project}, \cite[\href{https://stacks.math.columbia.edu/tag/0F6W}{Lemma 0F6W}]{stacks-project}, and Proposition \ref{prop:cotangent_rel_perf_animated}. 
In certain situations, some of these implications become equivalences. For instance:
\begin{itemize}
\item If $R \to S$ is a map of Noetherian rings, then the vanishing of the cotangent complex is equivalent to formal \'etaleness.
\item If $R \to S$ is a map of Noetherian \textit{F}-finite rings, then $S$ is relatively perfect over $R$ if and only if the cotangent complex is acyclic; see Theorem~\ref{thm:acyclic_cotanget_implies_rel_perf}.
\end{itemize}
\end{ex}

By a result of Gabber \citep[Proposition 5.2]{Kat86}, any relatively perfect scheme over a regular locally Noetherian scheme is flat. This does not hold in general if the base is not Noetherian.

\begin{ex}
\label{ex:rel_perf_not_flat_der}
Let $R$ be a discrete commutative $\FF_p$-algebra with flat Frobenius homomorphism, and let $I\subset R$ be an ideal satisfying $I^{[p]}=I$. Then $R/I$ is relatively perfect over $R$, but it is not necessarily flat. For example, take $R=\FF_p[x^{\frac{1}{p^{\infty}}}]$ with $I=(x^{\frac{1}{p^{\infty}}})$. In particular, a relatively perfect map need not be weakly \'etale.
\end{ex}

The proof of the following lemma is a direct application of Remark \ref{rem:basic_properties_general_der} and is therefore omitted.

\begin{lem}
\label{lem:relative_perfect_basic_properties_der}
\begin{enumerate}
\item[(i)] Let $R\to S \to T$ be maps of derived $\FF_p$-algebras. If $R\to S$ is relatively perfect, then $T$ is relatively perfect over $R$ if and only if it is relatively perfect over $S$.
\item[(ii)] Let $R$ be an derived $\FF_p$-algebra and let $R'$ and $S$ be animated commutative $R$-algebras. If $S$ is relatively perfect over $R$, then $S'=S\otimes^L_R R'$ is relatively perfect over $R'$.
\end{enumerate}
\end{lem}

The following two statements are direct analogues of Lemmas 3.16 and 3.18 in the work of Bhatt and Scholze \citep{BS17}.

\begin{prop}
\label{prop:rel_perf_tor-ind}
Let $R$ be a derived $\FF_p$-algebra with flat absolute Frobenius map. If both $S$ and $T$ are relatively perfect $R$-algebras, then $S\otimes_R^L T$ is coconnective. In particular, if $R$, $S$ and $T$ are discrete, then 
\[
S\otimes_R^L T \simeq S\otimes_R T.
\]
\end{prop}
\begin{proof}
By Lemma \ref{lem:relative_perfect_basic_properties_der}, the derived tensor product $S\otimes_R^LT$ is relatively perfect over $R$ and hence coconnective by Lemma \ref{lem:flat_frobenius_implies_coconnective}. If $R$, $S$ and $T$ are discrete, then the tensor product is both connective and coconnective, and is therefore discrete.
\end{proof}

\begin{cor}
\label{cor:perfect_base-change_der}
Let $Y$ be a discrete $\FF_p$-scheme with a flat Frobenius endomorphism, and consider the following pullback square of (discrete) $\FF_p$-schemes.
\[
\begin{tikzcd}
{X'} \arrow[r, "g'"] \arrow[d, "f'"] & {X} \arrow[d, "f"] \\
{Y'} \arrow[r, "g"]           & {Y}          
\end{tikzcd}
\]
If $X$ and $Y'$ are relatively perfect over $Y$ then, for any $K\in D_{\operatorname{qc}}(X)$, the base change morphism
\[
Lg^*Rf_*K\to Rf'_*L{g'}^*K
\]
is an isomorphism in $D_{\operatorname{qc}}(Y')$.
\end{cor}
\begin{proof}
Similarly as Lemma 3.18 is deduced from Lemma 3.16 in \citep{BS17}, the statement here reduces to the preceding Tor-independence result.
\end{proof}

The following result shows that the cotangent complex of a relatively perfect map of animated rings vanishes. We note that a similar result appears in \citep[Corollary 3.8]{Bha12}, where Bhatt proves the vanishing of the cotangent complex for relatively perfect maps in the setting of simplicial commutative rings. In \citep[Proposition 5.13]{Sch12} Scholze proves it for discrete rings.

\begin{prop}
\label{prop:cotangent_rel_perf_animated}
Let $R\to S$ be a map of animated commutative $\FF_p$-algebras. If $S$ is relatively perfect over $R$, then $L_{S/R}$ is acyclic.
\end{prop}
\begin{proof}
The cotangent complex of an equivalence is acyclic. Consequently, the proposition follows from Lemma \ref{lem:cotangent_relative_frobenius}.
\end{proof}


Lurie shows in \citep[Corollary 25.3.6.6]{Lur18} that a map of animated rings is an equivalence if and only if its cotangent complex vanishes and the map on $H^0$ is an isomorphism. Using this result, we deduce the following corollary to the previous proposition.

\begin{cor}
\label{cor:equivalence_animated}
A map $R\to S$ of animated commutative $\FF_p$-algebras is an equivalence if and only if it is relatively perfect and the induced map $H^0(R) \to H^0(S)$ is an isomorphism of commutative rings.
\end{cor}


The following example gives a different proof of \citep[Proposition 11.6]{BS17}.

\begin{ex}
Any map $R\to S$ of perfect animated $\FF_p$-algebras is relatively perfect, and thus is an equivalence if and only if it is an isomorphism on $H^0$. Applying this argument to the natural map $R\to H^0(R)$ shows that $R$ is discrete. 
\end{ex}

\begin{rem}
\label{rem:relative_perfect_over_regular}
Relatively perfect algebras over regular rings have at least two structural properties that, in the Noetherian case, are known to characterize regularity: Kunz’s theorem states that a Noetherian $\FF_p$-algebra $R$ is regular if and only if its absolute Frobenius map is flat \citep[Theorem 2.1]{Kun69}, and by a theorem of Andr\'e \citep[Supplement Theorem 30]{And74}, regularity of $R$ is equivalent to $L_{R/\FF_p}\simeq \Omega_{R/\FF_p}[0]$ with $\Omega_{R/\FF_p}$ being a flat $R$-module.\footnote{See also Lurie’s derived version of Popescu’s theorem \citep[Theorem 3.7.5]{Lur04} for a generalization to animated rings.}

As we have seen, any relatively perfect algebra over a discrete regular Noetherian base ring has both a flat Frobenius endomorphism and a flat cotangent complex — despite not being Noetherian in general. This indicates that relatively perfect algebras over regular rings could be seen as part of a class of rings that generalize certain properties of regularity beyond the Noetherian setting.

\end{rem}

\subsection*{Relative Inverse Limit Perfection} Our next goal is to construct an explicit right adjoint to the inclusion $\dalgperf{R}\into \dalg_R$ assuming that the functor $-\otimes_{R,F}R$ preserves inverse limits.

To streamline several upcoming statements, we introduce the following two notions. For a general treatment of dualizable objects, we refer to \citep[Sections 4.6.1 and 4.6.2]{Lur17}. We mainly rely on the fact that tensoring with dualizable modules preserves small limits.

\begin{defi}
\label{defi:F-dualizable-F-finite-projective}
A discrete commutative $\FF_p$-algebra $R$ is called \emph{\textit{F}-finite-projective} if and only if its Frobenius endomorphism makes its target $F_*R$ a finite projective module over its source $R$.

A derived $\FF_p$-algebra $R$ is called \emph{\textit{F}-dualizable} if and only if its Frobenius endomorphism makes its target $F_*R$ a dualizable module over its source $R$.
\end{defi}

Any relatively perfect algebra over a regular Noetherian \textit{F}-finite ring is \textit{F}-finite-projective. Moreover, the property of being \textit{F}-finite-projective is preserved under free extensions of finite type.

As shown by Lurie, dualizability of modules over animated rings coincides with being perfect and compact in the derived category. This yields the following reformulation of \textit{F}-dualizability in the connective case.

\begin{rem}
\label{rem:f-compact}
Let $R$ be an animated commutative $\FF_p$-algebra. By \citep[Proposition 6.2.6.2]{Lur18} and \citep[Proposition 7.2.4.2]{Lur17}, a derived $R$-module $M$ is dualizable if and only if it is perfect, which in turn is equivalent to being a compact object of the $\infty$-category of derived $R$-modules. In particular, any discrete \textit{F}-finite-projective ring is \textit{F}-dualizable.
\end{rem}

\begin{lem}
\label{lem:co_limit_of_perfect}
For a derived $\mathbb{F}_p$-algebra $R$, the category $\dalgperf{R}$ is closed under small colimits in $\dalg_R$. Moreover, if $R$ is \textit{F}-dualizable, it is also closed under small limits.
\end{lem}
\begin{proof}
We show both statements using the description of the relative Frobenius maps given in Remark \ref{rem:limit-colimit-frobenius_der}. The case of colimits is immediate. For limits, let $S:I\to \dalgperf{R}$ be a diagram of relatively perfect $R$-algebras. To show that the limit is relatively perfect, it suffices to prove that the natural map 
$(\lim_{i\in I} S_i) \otimes_{R,F}R \to \lim_{i\in I} (S_i \otimes_{R,F}R)$ is an equivalence. This reduces to showing that the underlying map of derived $R$-modules is an equivalence. Since $F_*R$ is a dualizable $R$-module, the functor $(-)\otimes_{R,F}R$ preserves small limits of derived $R$-modules.
\end{proof}

\begin{lem}
\label{lem:presentable}
For a derived $\mathbb{F}_p$-algebra $R$, the category $\dalgperf{R}$ is presentable.
\end{lem}
\begin{proof}
Composing the absolute Frobenius functor constructed in \citep[Construction~2.4.1]{Hol23} with the evaluation at the loop $1:*\to *$ in $B\NN$ yields a functor
\[
\dalg_{\FF_p} 
\to 
\fun(B\NN, \dalg_{\FF_p}) 
\to 
\fun(\Delta^1, \dalg_{\FF_p}),
\]
which assigns to a derived $\FF_p$-algebra its absolute Frobenius map. This induces a composite functor
\[
\dalg_R
\to 
\fun(\Delta^1, \dalg_{\FF_p}) 
\to
\fun(\Delta^1, \fun(\Delta^1, \dalg_{\FF_p}))
\simeq
\fun(\Delta^1\times \Delta^1, \dalg_{\FF_p}) 
\]
which sends an $R$-algebra $S$ to the naturality square of the Frobenius maps of $R\to S$.
We obtain a cartesian diagram of $\infty$-categories
\[
\begin{tikzcd}
\dalgperf{R} \arrow[r, hook] \arrow[d]              & \dalg_R \arrow[d]                               \\
{\fun(\Lambda^2_0, \dalg_{\FF_p}) } \arrow[r, hook] & {\fun(\Delta^1\times \Delta^1, \dalg_{\FF_p}) }
\end{tikzcd}
\]
where the lower horizontal arrow is fully faithful with essential image the full subcategory of cocartesian squares in $\dalg_{\FF_p}$. By \citep[Remark~4.2.23]{Rak20}, and  \citep[Proposition~5.5.3.6 and Proposition~5.5.3.11]{Lur09}, all categories in the diagram are presentable, except possibly $\dalgperf{R}$. The lower horizontal and right vertical functors preserve filtered colimits and are therefore accessible. Consequently, by \citep[Proposition~5.4.6.6]{Lur09}, the category $\dalgperf{R}$ is accessible and thus presentable by Lemma \ref{lem:co_limit_of_perfect}.
\end{proof}

\begin{prop}
\label{prop:adjoint_functor_thm}
For any derived $\FF_p$-algebra $R$, the inclusion 
\[
\dalgperf{R}\into \dalg_R
\]
admits a right adjoint. Moreover, if $R$ is $F$-dualizable, then the inclusion admits both a left and a right adjoint.
\end{prop}
\begin{proof}
This follows from Lemmas~\ref{lem:co_limit_of_perfect} and~\ref{lem:presentable}, and from the adjoint functor theorem \citep[Corollary~5.5.2.9]{Lur09}.
\end{proof}

We next construct an explicit right adjoint provided that the base is $F$-dualizable.

\begin{con}
\label{con:derived_rel_frobenius_tower}
For a map $R\to S$ of derived $\FF_p$-algebras, the following commutative diagram
\[
\begin{tikzcd}
                & \dots \arrow[r]                              & R \arrow[r, "\id_R"]                                  & R \arrow[r, "\id_R"]                                & R \arrow[r, "\id_R"]                     & R \\
\dots \arrow[r] & R \arrow[r, "F"] \arrow[ru, "F^3"] \arrow[d] & R \arrow[r, "F"] \arrow[ru, "F^2"] \arrow[d] & R \arrow[r, "F"] \arrow[ru, "F"] \arrow[d] & R \arrow[d] \arrow[ru, "\id"] &   \\
\dots \arrow[r] & S \arrow[r, "F"]                             & S \arrow[r, "F"]                             & S \arrow[r, "F"]                           & S                               &  
\end{tikzcd}
\]
induces an inverse system of derived $R$-algebras:
\[
\dots\to S\otimes^L_{R,F^3}R\to[F\otimes \id_R] S\otimes^L_{R,F^2}R\to[F\otimes \id_R] S\otimes^L_{R,F}R \to[F\otimes \id_R] S\otimes^L_{R}R\simeq S.
\]
For each $n\ge0$, the map $S\otimes^L_{R,F^{n+1}}R\to[F\otimes \id_R] S\otimes^L_{R,F^n}R$ is 
a relative Frobenius map associated with the second factor map $R \to S\otimes_{R,F^n}R$.
\end{con}

\begin{defi}
\label{def:derived_rel_frobenius_tower}
Let $R\to S$ be a map of derived $\FF_p$-algebras. We call the inverse system of derived $R$-algebras 
\[
\dots \to[F\otimes\id_R] S\otimes^L_{R,F^2}R\to[F\otimes \id_R] S\otimes^L_{R,F}R \to[F\otimes \id_R] S\otimes^L_{R}R\simeq S
\]
the \emph{relative Frobenius tower of $S$ over $R$}. We denote by $S^{\perfection/R}$ the inverse limit over this system in the category of derived $R$-algebras.
\end{defi}

For an animated ring $R$ and elements $f_1, \dots, f_n \in H^0(R)$, we write
\[
R\sslash(f_1, \dots, f_n) = \FF_p \otimes^L_{\FF_p[x_1, \dots, x_n]} R,
\]
where the map $ \FF_p[x_1, \dots, x_n] \to R $ corresponds to the $f_i$, and $ \FF_p[x_1, \dots, x_n] \to \FF_p $ sends all variables to zero.

\begin{ex}[Derived Completion]
\label{ex:derived-completion}
Let $R$ be an animated $\FF_p$-algebra, and let $I=(f_1,\dots, f_n)\subset H^0(R)$ be a finitely generated ideal for some $n\in \NN$. Then there are canonical equivalences
\begin{align*}
R^{\wedge}_I 
&\simeq \varprojlim_k \left(\FF_p[x_1,\dots x_n]/(x_1,\dots,x_n)^{[p^k]}\otimes^L_{\FF_p[x_1,\dots x_n]} R\right)\\
&\simeq \varprojlim_k \left(\FF_p\otimes^L_{\FF_p[x_1,\dots x_n], F^k}\FF_p[x_1,\dots x_n] \otimes_{\FF_p[x_1,\dots x_n]}^L R\right)\\
&\simeq \varprojlim_k \left((\FF_p\otimes^L_{\FF_p[x_1,\dots x_n]}R) \otimes_{R, F^k}^L R\right)\\
&\simeq (R\sslash(f_1,\dots,f_n))^{\perfection/R},
\end{align*}
showing that the derived completion of an animated $\FF_p$-algebra can be realized as the limit of a relative Frobenius tower.
\end{ex}

The limit of the relative Frobenius tower depends on the base only up to a relatively perfect map.

\begin{lem}
\label{lem:independence_of_base}
Let $R\to T\to S$ be maps of derived $\FF_p$-algebras. If $T$ is relatively perfect over $R$, then the canonical map is a equivalence $S^{\perfection/R}\simeq S^{\perfection/T}$ of derived $R$-algebras.
\end{lem}
\begin{proof}
The assumption implies that the relative Frobenius towers over $R$ and $T$ are equivalent, so their limits are equivalent as well.
\end{proof}

If the base has a flat absolute Frobenius map, then the limit over the relative Frobenius tower is coconnective. The additional assumption in the following lemma is satisfied, for instance, when $R \to S$ is a relatively semiperfect map of animated $\mathbb{F}_p$-algebras; see Definition \ref{def:rel_semiperf_der}.

\begin{lem}
\label{lem:perfection_sees_only_H0_flat}
Let $R\to S$ be a map of derived $\FF_p$-algebras. If the absolute Frobenius map of $R$ is flat, then $S^{\perfection/R}$ is coconnective. Moreover, if the relative Frobenius map of $H^0(R)\to H^0(S)$ is surjective, then there is a canonical isomorphism $H^0(S^{\perfection/R})\cong (H^0S)^{\perfection/H^0R}$.
\end{lem}
\begin{proof}
Let $k \le -1$. By Lemma \ref{lem:rel_frob_zero_on_cohomology}, applying $H^k$ to the relative Frobenius tower yields a system whose transition maps are zero:
\[
\dots \to H^k(S)\otimes_{H^0(R),F^2}H^0(R) \to H^k(S)\otimes_{H^0(R),F}H^0(R) \to H^k(S)
\]
Hence, the limit of this system vanishes in $D(R)$, which implies that both outer terms in the Milnor short exact sequence vanish:
\[
0
\to 
{\varprojlim_n}^1 (H^{k-1}S\otimes_{H^0(R),F^n}H^0R)
\to H^k(S^{\perfection/R}) 
\to 
{\varprojlim_n}^0 (H^k S\otimes_{H^0(R),F^n}H^0R) 
\to 
0.
\]
For $i = 0, 1$, we denote by $ {\varprojlim}_n^i $ the $i$-th cohomology of the limit of the corresponding diagram in the $\infty$-category $ D(\FF_p) $.

For $H^0$, the Milnor exact sequence yields an isomorphism: 
\[
H^0(S^{\perfection/R}) 
\cong
{\varprojlim_n}^0 (H^0(S)\otimes_{H^0(R),F^n}H^0(R)). 
\]
If the relative Frobenius map of $H^0(R)\to H^0(S)$ is surjective, then the transition maps in the inverse system
\[
\dots \to[F\otimes\id ] H^0(S)\otimes^L_{H^0(R),F^2}H^0(R)\to[F\otimes \id] H^0(S) \otimes^L_{H^0(R),F}H^0(R) \to H^0(S)
\]
are surjective (see the proof of Lemma \ref{lem:connective_surjective_der}) and thus the tower satisfies the Mittag-Leffler condition. This gives a canonical isomorphism:
\[
{\varprojlim_n}^0 (H^0(S)\otimes_{H^0(R),F^n} H^0(R)) \cong (H^0 S)^{\perfection/H^0(R)}.
\]
\end{proof}

If the base is \textit{F}-dualizable, then the limit over the relative Frobenius tower is relatively perfect.

\begin{lem}
\label{lem:derived_relative_perfection}
Let $R$ be an \textit{F}-dualizable derived $\FF_p$-algebra. Then, for any derived $R$-algebra $S$, the derived $R$-algebra $S^{\perfection/R}$ is relatively perfect.
\end{lem}
\begin{proof}
As noted in Remark \ref{rem:limit-colimit-frobenius_der}, the composition
\[
S^{\perfection/R}\otimes_{R,F}R 
\to 
\varprojlim\left(\dots\to[F\otimes \id_R] S\otimes^L_{R,F^2}R\to[F\otimes \id_R] S\otimes^L_{R,F}R\right) 
\to 
S^{\perfection/R}
\]
is a relative Frobenius map for $R\to S^{\perfection/R}$. The first map is an equivalence due to the dualizability assumption on $R$, while the second map is an equivalence because the limit of the tower does not depend on finitely many initial terms.
\end{proof}

Under the assumptions of the previous lemma, the construction $S \mapsto S^{\perfection/R}$ induces a functor 
\[
(-)^{\perfection/R}:\dalg_R \to \dalgperf{R}
\] 
which we refer to as the \emph{relative inverse limit perfection over $R$}.

\begin{prop}
\label{prop:derived_relative_inverse_right_adjoint}
Let $R$ be an \textit{F}-dualizable derived $\FF_p$-algebra. Then, the relative inverse limit perfection $(-)^{\perfection/R}$ is right adjoint to the inclusion $\dalgperf{R}\subset \dalg_R$.
\end{prop}
\begin{proof}
Let $\varepsilon: (-)^{\perfection/R}\to \id_{\dalg_R}$ be the natural transformation arising from the limit structure. We show that for $T\in \dalgperf{R}$ and $S\in \dalg_R$, the composition
\[
\operatorname{Map}_{\dalgperf{R}}(T,S^{\perfection/R})\to[\simeq] \operatorname{Map}_{\dalg_R}( T, S^{\perfection/R}) \to[\varepsilon_S\circ] \operatorname{Map}_{\dalg_R}(T, S)
\]
is an equivalence of anima, which means that $\varepsilon$ serves as the counit of an adjunction between the inclusion and $(-)^{\perfection/R}$. Since $\dalgperf{R}$ is a full subcategory of $\dalg_R$, it suffices to show that the map
\[
\varprojlim_n \operatorname{Map}_{\dalg_R}(T, S\otimes_{R,F^n}R) \simeq \operatorname{Map}_{\dalg_R}(T, S^{\perfection/R}) \to \operatorname{Map}_{\dalg_R}(T, S)
\]
is an equivalence. This follows from the fact that $T$ is relatively perfect over $R$, ensuring that each map in the inverse system of anima is an equivalence; see Remark \ref{rem:equivalence_of_spaces_der}.
\end{proof}

\begin{rem}
The relative perfection constructed by Kato in \citep{Kat86} corresponds to a relative colimit perfection in the category of commutative algebras and serves as a left adjoint to the inclusion functor when the base ring is \textit{F}-finite-projective.
\end{rem}

\section{Relative Semiperfectness and F-finiteness}
\label{sec:rel_semiperf}
We now consider connective rings. Let $R\to S$ be a map of animated $\FF_p$-algebras. Since the forgetful functor $\dalg_R\to D(R)$ preserves small limits, the $k$-th cohomology group of $S^{\perfection/R}$ fits into a Milnor exact sequence for each $k\in \ZZ$:
\[
0
\to 
{\varprojlim_n}^1 H^{k-1}(S\otimes^L_{R,F^n}R)
\to H^k(S^{\perfection/R}) 
\to 
{\varprojlim_n}^0 H^k(S\otimes^L_{R,F^n}R) 
\to 
0.
\]
As each level of the relative Frobenius tower is connective, it follows that $S^{\perfection/R}$ is always cohomologically $1$-connective. We now explore a condition ensuring that $S^{\perfection/R}$ remains connective.

We recall that an animated $\FF_p$-algebra is called \emph{semiperfect} if and only if its Frobenius endomorphism is surjective on $H^0$.

\begin{defi}
\label{def:rel_semiperf_der}
A map $R\to S$ of animated $\FF_p$-algebras is called \emph{relatively semiperfect} if and only if its relative Frobenius map $F_{S/R}$ is surjective on $H^0$. In this case, we call $S$ a \emph{relatively semiperfect $R$-algebra}.

We denote by $\sperf{R}\subset \ACR_R$ the full subcategory of relatively semiperfect $R$-algebras.
\end{defi}

\begin{ex}
Any map $R\to S$ that is surjective on $H^0$ is relatively semiperfect. 
\end{ex}

Relative semiperfectness ensures not only the connectivity of the relative inverse limit perfection but also that of the fiber of the natural map
$\varepsilon: (-)^{\perfection/R}\to \id$.

\begin{lem}
\label{lem:connective_surjective_der}
Let $R$ be an animated $\FF_p$-algebra. Then, for any relatively semiperfect $R$-algebra $S$, the following statements hold:
\begin{enumerate}
\item[(i)] The derived $R$-algebra $S^{\perfection/R}$ is connective; that is, an animated ring.
\item[(ii)] The natural map $\varepsilon_S: S^{\perfection/R}\to S$ is surjective on $H^0$.
\end{enumerate}
\end{lem}
\begin{proof}
The map $S\otimes^L_{R,F^{n+1}}R\to[F\otimes \id_R] S\otimes^L_{R,F^n}R$ in the relative Frobenius tower is a relative Frobenius map for $R\to S\otimes_{R,F^n}R$. Hence, part (ii) of Remark \ref{rem:basic_properties_general_der}, applied to $R\ot[F^n]R\to S$, shows that
\[
\begin{tikzcd}
{S\otimes_{R,F}R} \arrow[r] \arrow[d, "F_{S/R}"'] & {S\otimes_{R,F^{n+1}}R} \arrow[d, "F\otimes \id"] \\
S \arrow[r]                                       & {S\otimes_{R,F^{n}}R}                            
\end{tikzcd}
\]
is a pushout square of animated rings. In particular, if $R\to S$ is relatively semiperfect, each map in the inverse system is surjective on $H^0$.
Now, both statements follow from Milnor exact sequences. For (i), since each stage of the relative Frobenius tower is connective, the inverse limit is always cohomologically $1$-connective. Moreover, as argued above, the induced system on $H^0$ satisfies the Mittag-Leffler condition. Consequently, the left-hand and right-hand sides of the short exact sequence vanish:
\[
0\to {\varprojlim_n}^1 H^0(S\otimes_{R,F^n}R) 
\to 
H^1(S^{\perfection/R}) 
\to 
{\varprojlim_n}^0 H^1(S\otimes_{R,F^n}R)\to 0
\]
Regarding (ii), the canonical map on $H^0$ factors
 as
\[
H^0(S^{\perfection/R}) 
\to 
\varprojlim_n H^0(S\otimes_{R,F^n}R)
\to
H^0(S),
\]
and both homomorphisms are surjective.
\end{proof}

For an animated $\FF_p$-algebra $R$ we denote by $\perf{R}\subset \ACR_R$ the full subcategory consisting of relatively perfect animated $R$-algebras. As in Lemmas~\ref{lem:co_limit_of_perfect} and~\ref{lem:presentable}, one verifies that $\perf{R}$ is presentable, that the inclusion always preserves small colimits, and that it preserves small limits if $R$ is \textit{F}-dualizable. Hence, by Proposition~\ref{prop:adjoint_functor_thm}, the inclusion admits a right adjoint, and under the additional assumption that $R$ is \textit{F}-dualizable, also a left adjoint. If the base is $F$-dualizable, the relative inverse limit perfection is a right adjoint to the inclusion $\perf{R}\subset \sperf{R}$.

\begin{prop}
\label{prop:right_adjoint_animated}
Let $R$ be an $F$-dualizable animated $\FF_p$-algebra. Then, the relative inverse limit perfection induces a functor 
\[
(-)^{\perfection/R}: \sperf{R} \to \perf{R}
\]
which is right adjoint to the inclusion $\perf{R}\subset \sperf{R}$.
\end{prop}
\begin{proof}
Let $S$ be a relatively semiperfect $R$-algebra. By Lemma \ref{lem:derived_relative_perfection}, the derived $R$-algebra $S^{\perfection/R}$ is relatively perfect over $R$, and by Lemma \ref{lem:connective_surjective_der}, it is connective. The proof now follows similarly to that of Proposition \ref{prop:derived_relative_inverse_right_adjoint}, as all involved categories are full subcategories of $\dalg_R$. Once again, we use the natural transformation $\varepsilon$ from the limit structure as the counit.
\end{proof}

The following results compare the relative inverse limit perfection with some notions of completion.

\begin{lem}
\label{lem:derived-compl-flat-frob}
Let $R$ be an $\FF_p$-algebra with a flat Frobenius endomorphism, and let $I\subset R$ be a finitely generated ideal. Then, the derived completion $R^{\wedge}_I$ canonically identifies with the classical $I$-adic completion of $R$.
\end{lem}
\begin{proof}
Let $f_1,\dots,f_n\in R$ elements that generate $I$. We have the following canonical equivalences:
\[
R^{\wedge}_I \simeq (R\sslash(f_1,\dots,f_n))^{\perfection/R} \simeq (R/I)^{\perfection/R}.
\]
Indeed, the first equivalence holds by Example \ref{ex:derived-completion}, and the second equivalence holds by Lemmas \ref{lem:perfection_sees_only_H0_flat} and \ref{lem:connective_surjective_der} as $R$ has a flat absolute Frobenius map. Now, the right-hand-side is given by $\varprojlim_k(R/I^{[p^k]})$ which naturally identifies with the classical $I$-adic completion of $R$, as $I$ is finitely generated.
\end{proof}

\begin{prop}
Let $R \to S$ be a map of animated $\FF_p$-algebras, and assume that the Frobenius endomorphism of $R$ is flat. If $R\to H^0(S)$ is surjective with finitely generated kernel $I$, then $S^{\perfection/R}$ is canonically equivalent to the Adams completion $\operatorname{Comp}(R\to S)$.
\end{prop}
\begin{proof}
By \citep[Proposition 3.2.5]{BBST}, the Adams completion $\operatorname{Comp}(R\to S)$ canonically identifies with the derived $I$-completion of $R$ and thus by Lemma \ref{lem:derived-compl-flat-frob} with the classical $I$-adic completion of $R$. On the other hand, Lemma \ref{lem:perfection_sees_only_H0_flat} shows $S^{\perfection/R}\simeq (H^0S)^{\perfection/R}$ which is also isomorphic to the classical completion.
\end{proof}

\begin{ex}
\label{ex:perfection_complete_Frobenius_powers_der}
(i) Let $R \to S$ be a map of animated $\FF_p$-algebras, and assume that the Frobenius endomorphism of $R$ is flat. Then by Lemma \ref{lem:perfection_sees_only_H0_flat}, we have a canonical equivalence: 
\[
S^{\perfection/R}\cong (H^0S)^{\perfection/R}.
\]
Let $I$ denote the kernel of the map $R\to H^0(S)$. Since this map factors as $R\to R/I \to H^0(S)$, the relative Frobenius tower of $H^0(S)$ over $R$ is given by the following:
\[
\dots
\to
H^0S\otimes_{R/I,F^3}^L R/I^{[p^3]}
\to[F\otimes \operatorname{\, pr}]
H^0S\otimes_{R/I,F^2}^L R/I^{[p^2]}
\to[F\otimes \operatorname{\, pr}]
H^0S\otimes_{R/I,F}^L R/I^{[p]}
\to
H^0S.
\]
In particular, if $R\to H^0(S)$ is surjective, then $S^{\perfection/R}\cong \varprojlim_n R/I^{[p^n]}$.

(ii) Let $R$ be an \textit{F}-finite-projective $\FF_p$-algebra, and let $S$ be a relatively semiperfect $R$-algebra. Then $S^{\perfection/R}$ is relatively perfect over $R$ and thus also \textit{F}-finite-projective. Let $I$ denote the kernel of the natural map $S^{\perfection/R}\to H^0(S)$. By Lemmas \ref{lem:perfection_sees_only_H0_flat} and \ref{lem:independence_of_base}, together with part (i), we obtain isomorphisms
\[
S^{\perfection/R}
\cong 
(H^0S)^{\perfection/R} 
\cong 
(H^0S)^{\perfection/S^{\perfection/R}} 
\cong 
\varprojlim_n S^{\perfection/R}/I^{[p^n]}
\]
showing that $S^{\perfection/R}$ is complete with respect to the topology induced by $\{I^{[p^n]}\}_{n\ge0}$ as a fundamental system of neighborhoods.
\end{ex}

The assumptions of the following proposition are satisfied, for instance, when $R$ is a regular Noetherian \textit{F}-finite ring and $H^0(S)$ is Noetherian; see Proposition~\ref{prop:perfection_is_Noetherian_der}.

\begin{prop}
\label{prop:adams_complete_der}
Let $R$ be an \textit{F}-finite-projective ring, and let $S$ be a relatively semiperfect animated $R$-algebra. If the kernel of the map
\[
S^{\perfection/R}\to S
\]
is finitely generated, then this map is Adams-complete.
\end{prop}
\begin{proof}
Consider the following maps $S^{\perfection/R} \to S \to H^0(S)$. By \citep[Corollary 3.2.10]{BBST}, the second map is universally Adams complete. Hence, by \citep[Proposition 3.2.8]{BBST}, it suffices to show that the composite map $S^{\perfection/R} \to H^0(S)$ is Adams complete. Now, by Lemma \ref{lem:perfection_sees_only_H0_flat}, we may assume that $S$ is discrete. In this case, Example \ref{ex:perfection_complete_Frobenius_powers_der} shows that $S^{\perfection/R}$ is complete with respect to the topology induced by $\{I^{[p^n]}\}_{n\ge0}$, where $I$ is the kernel of the map $S^{\perfection/R} \to S$. Since the ideal $I$ is finitely generated, the result follows from \citep[Proposition 3.2.5]{BBST}.
\end{proof}

\subsection{\textit{F-finiteness}}
We next consider \textit{F}-finite maps of animated $\FF_p$-algebras and examine their relation to relatively semiperfect maps.

\begin{defi}
\label{def:rel_f-finite_der}
A map $R\to S$ of animated $\FF_p$-algebras is called \emph{\textit{F}-finite} if the relative Frobenius map $F_{S/R}$ is finite on $H^0$. In this case, we also say that \emph{$S$ is a relatively \textit{F}-finite $R$-algebra}. When $R = \FF_p$, we simply say that \emph{$S$ is \textit{F}-finite}.
\end{defi}

\begin{ex}
If $R\to S$ is of finite type on $H^0$ then $S$ is relatively \textit{F}-finite over $R$.
\end{ex}

The next lemma follows directly from Remark \ref{rem:basic_properties_general_der}. It is also closely related to Lemma 2 in Hashimoto’s paper \citep{Has15}, as our notions of semiperfectness and \textit{F}-finiteness are defined entirely in terms of $H^0$.

\begin{lem}
\label{lem:basic_properties_P_der}
Let $P\in\{$relatively semiperfect, \textit{F}-finite$\}$. For maps of animated $\FF_p$-algebras, the following statements hold:
\begin{enumerate}
\item If both $A\to B$ and $A\to C$ satisfy $P$, then so does the composition. Moreover, if $A \to C$ satisfies $P$, then so does $B \to C$.
\item If $A \to B$ satisfies $P$, then for any map $A \to A^\prime$, the base change $A^\prime \to A^\prime \otimes^L_A B$ also satisfies $P$.
\item If both $A\to B$ and $A\to C$ satisfy $P$, then $A\to B\otimes^L_{A}C$ also satisfies $P$.
\end{enumerate}
\end{lem}

\begin{ex}
\label{ex:maps_of_semiperf_and_f-fin_der}
If $R \to S$ is a map of animated $\mathbb{F}_p$-algebras and $S$ is semiperfect (resp. \textit{F}-finite), then the map is relatively semiperfect (resp. relatively \textit{F}-finite).
\end{ex}

\begin{ex}
Let $R\to S$ be a map of animated $\FF_p$-algebras. If $S$ is relatively \textit{F}-finite over $R$, then any $S$-algebra of finite type is relatively \textit{F}-finite over $R$. This follows from Lemma \ref{lem:basic_properties_P_der} and the fact that $\FF_p[x]$ is \textit{F}-finite. 
\end{ex}

The following lemma relates relative semiperfectness to the vanishing of the degree-zero cohomology of the cotangent complex.

\begin{lem}
\label{lem:cotangent_semiperfect_der}
Let $R\to S$ be a map of animated $\FF_p$-algebras. If $S$ is relatively semiperfect over $R$, then $L_{S/R}$ is cohomologically $-1$-connective. 

Moreover, if $R\to S$ is \textit{F}-finite, then $S$ is relatively semiperfect over $R$ if and only if the module of Kähler differentials $\Omega_{H^0S/H^0R}\cong H^0(L_{S/R})$ vanishes.
\end{lem}
\begin{proof}
If $F_{S/R}$ is surjective on $H^0$, then $L_{F_{S/R}}$ is cohomologically $-1$-connective by \citep[Corollary 25.3.6.4]{Lur18}. Therefore,
\[
H^0(L_{S/R})=H^0(L_{F_{S/R}})=0
\]
by Lemma \ref{lem:cotangent_relative_frobenius}.

Now, suppose that $R\to S$ is \textit{F}-finite and that $L_{S/R}$ is $-1$-connective. Let $A= H^0(R)$ and $B=H^0(S)$. We have to show $B^p[A]=B$. Since 
\[
\Omega_{B/B^p[A]} \cong \Omega_{B/A} = 0,
\]
this follows from \citep[Corollaire 21.1.6]{EGAIV}.
\end{proof}

Analogous to a finite type map, an \textit{F}-finite map admits a factorization into a finite type free followed by an “\textit{F}-surjective” map. 

\begin{lem}
\label{lem:f-finite_smooth_f-surjective}
Let $R\to S$ be an \textit{F}-finite map of animated $\FF_p$-algebras. Then there exists an integer $n \ge 0$ and a factorization
\[
R \to R[x_1,\dots, x_n] \to S,
\]
where the first map is free of finite type and the second is relatively semiperfect.
\end{lem}
\begin{proof}
Let $A=H^0(R)$ and $B=H^0(S)$. If $S$ is relatively \textit{F}-finite over $R$, then there exist finitely many elements $s_1, \dots, s_n \in B$ that generate $B$ as a $B^p[A]$-algebra. These elements define a map from a free $R$-algebra $R[x_1,\dots, x_n]$ to $S$, which is relatively semiperfect since 
\[
B^p[H^0(R[x_1,\dots,x_n])]=B^p[A,s_1,\dots,s_n]=H^0(S).
\]
\end{proof}

If the base field in the following proposition is perfect, then the result is due to Kunz; see \citep[Proposition 1.1]{Kun76}.

\begin{prop}
Let $k$ be field of characteristic $p$, and let $R$ be a discrete Noetherian $k$-algebra such that the field extensions $k \to R/m$ are separable for all maximal ideals $m \subset R$. If $R$ is relatively $F$-finite over $k$, then $R$ has finite Krull dimension. More precisely, if $k[x_1,\dots,x_n] \to R$ is a relatively semiperfect map of $k$-algebras for some $n \ge 0$, then the Krull dimension of $R$ is at most~$n$.
\end{prop}
\begin{proof}
Since localizations are relatively perfect, we may assume that $(R,m)$ is a local ring such that the field extension $k\into K=R/m$ is separable. The canonical fiber sequence $L_{R/k}\otimes^L_R K \to L_{K/k} \to L_{K/R}$ induces an exact sequence of $K$-vector spaces
\[
H^{-1}(L_{K/k})\to H^{-1}(L_{K/R})\to \Omega_{R/k}\otimes_R K \to \Omega_{K/k} \to \Omega_{K/R} \to 0.
\]
Since $k\into K$ is separable—or equivalently, formally smooth by \cite[\href{https://stacks.math.columbia.edu/tag/0321}{Lemma 0321}]{stacks-project}—we have $H^{-1}(L_{K/k})=0$ by \cite[\href{https://stacks.math.columbia.edu/tag/031J}{Proposition 031J}]{stacks-project}. Furthermore, as $R\to K$ is surjective, we have $\Omega_{K/R}=0$ and $H^{-1}(L_{K/R})=m/m^2$. Hence, we have a short exact sequence
\[
0\to m/m^2 \to \Omega_{R/k} \otimes_R K \to \Omega_{K/k} \to 0
\]
which shows $\dim(R) \le \dim_K (m/m^2) \le \dim_K (\Omega_{R/k} \otimes_R K)$. Moreover, since $k[x_1,\dots, x_n]\to R$ is relatively semiperfect, the cotangent complex $L_{R/k[x_1,\dots, x_n]}$ is cohomologically $-1$-connective. Hence, the long exact cohomology sequence associated to
\[
L_{k[x_1,\dots, x_n]/k}\otimes^L_{k[x_1,\dots, x_n]}R
\to 
L_{R/k}
\to
L_{R/k[x_1,\dots, x_n]}
\]
shows that $R^n\cong \Omega_{k[x_1,\dots, x_n]/k}\otimes_{k[x_1,\dots, x_n]}R \to \Omega_{R/k}$ is surjective. Consequently, we have
\[
\dim_K (\Omega_{R/k} \otimes_R K)\le n.
\]
\end{proof}

\begin{ex}
\label{ex:gabbers_construction_der}
Let $S$ be a discrete \textit{F}-finite $\FF_p$-algebra, and let  $R=\FF_p[x_1,\dots,x_n]\to S$ be a relatively semiperfect map. Denote by $s_i$ the image of $x_i$ under this map. Then, for all $k\ge0$, we obtain commutative diagrams 
\[
\begin{tikzcd}
{S[x_1,\dots,x_n]/(x_1^{p^{k+1}}-s_1,\dots,x_n^{p^{k+1}}-s_n)} \arrow[r] \arrow[d, "\cong"] & {S[x_1,\dots,x_n]/(x_1^{p^{k}}-s_1,\dots,x_n^{p^{k}}-s_n)} \arrow[d, "\cong"] \\
{S\otimes_{R,F^{k+1}}R} \arrow[r, "F\otimes \id"]           & {S\otimes_{R,F^{k}}R}          
\end{tikzcd}
\]
where the upper horizontal map is the Frobenius on $R$ and sends $x_i$ to $x_i$. Since $R$ is regular, the tensor products in the lower row are automatically derived.


Thus, in this case, the relative Frobenius tower coincides exactly with the inverse system constructed in \citep[Remark 13.6]{Gab04}. In particular, if $S$ is Noetherian, then $S^{\perfection/R}$ is a regular Noetherian ring by loc.~cit.\footnote{See also \citep[Theorem 10.9]{MP20} and \citep[Construction 2.2.3]{BBST} for discussions of Gabber’s result in the book by Ma–Polstra and in the article by Bhatt–Blickle–Schwede–Tucker, respectively.}
\end{ex}

\begin{prop}
\label{prop:perfection_is_Noetherian_der}
Let $R$ be a regular Noetherian \textit{F}-finite $\FF_p$-algebra, and let $R\to S$ be a relatively semiperfect map to an \textit{F}-finite animated $\FF_p$-algebra. If $H^0(S)$ is Noetherian, then $S^{\perfection/R}$ is regular Noetherian.
\end{prop}
\begin{proof}
By Lemma~\ref{lem:perfection_sees_only_H0_flat}, we may assume that $S$ is discrete. By Kunz’s theorem \citep[Theorem 2.1]{Kun69}, it suffices to show that $S^{\perfection/R}$ is Noetherian.

Let $n\ge0$, and let $A=\FF_p[x_1,\dots,x_n]\to R$ be a relatively semiperfect map. By Example~\ref{ex:gabbers_construction_der} and \citep[Remark 13.6]{Gab04}, the relatively perfect $A$-algebras $R^{\perfection/A}$ and $S^{\perfection/A}$ are regular Noetherian.

The map $R^{\perfection/A}\to S^{\perfection/A}$ is relatively perfect by Lemma \ref{lem:relative_perfect_basic_properties_der}, hence flat by \citep[Proposition 5.2]{Kat86}. Moreover, by Lemma \ref{lem:connective_surjective_der}, the map $R^{\perfection/A}\to R$ is surjective. Hence,
\[
S^{\perfection/A}\otimes_{R^{\perfection/A}} R\simeq S^{\perfection/A}\otimes^L_{R^{\perfection/A}} R
\]
is a discrete Noetherian $\FF_p$-algebra. Furthermore, as $R^{\perfection/A}$ is regular, $R\in D(R^{\perfection/A})$ is a perfect complex, so $-\otimes^L_{R^{\perfection/A}} R$ preserves limits. This yields the following canonical equivalences:
\begin{align*}
S^{\perfection/A}\otimes^L_{R^{\perfection/A}} R 
&\simeq 
S^{\perfection/R^{\perfection/A}}\otimes^L_{R^{\perfection/A}} R
\simeq 
(S\otimes^L_{R} R \otimes^L_{R^{\perfection/A}}R)^{\perfection/R} \\
&\simeq 
(H^0(S\otimes^L_{R} R \otimes^L_{R^{\perfection/A}}R))^{\perfection/R}
\simeq 
S^{\perfection/R}.
\end{align*}
Here, the first equivalence follows from Lemma \ref{lem:independence_of_base}, the second from the discussion above, and the third from Lemma \ref{lem:perfection_sees_only_H0_flat}.
\end{proof}

The following proposition can also be deduced from  \citep[Theorem 6.6.4.1]{Lur18} and \citep[Remark 13.6]{Gab04}.

\begin{prop}
\label{prop:dual_complex}
Noetherian\footnote{An animated ring $R$ is called \emph{Noetherian} if $H^0(R)$ has the same property and $H^k(R)$ is a finitely generated $H^0(R)$ module for all $k\le 0$.} \textit{F}-finite animated rings admit dualizing complexes.
\end{prop}
\begin{proof}
Let $S$ be a Noetherian \textit{F}-finite animated ring, and choose a relatively semiperfect map $R=\FF_p[x_1,\dots,x_n] \to S$. By Proposition \ref{prop:perfection_is_Noetherian_der}, $S^{\perfection/R}$ is regular Noetherian. The cohomology groups of $S$ are finitely generated over $H^0(S)$ and therefore also over $S^{\perfection/R}$. This shows that $S$ is an almost-perfect $S^{\perfection/R}$-module; see \citep[Proposition 7.2.4.17]{Lur17}. Now, if $\omega_{S^{\perfection/R}}^{\bullet}$ is a dualizing complex over $S^{\perfection/R}$, then $\Hom_T(S,\omega_T^{\bullet})$ is a dualizing complex over $S$ by \citep[Corollary 6.6.3.3]{Lur18}.
\end{proof}

\section{Factoring Maps of F-finite Animated Rings}
\label{sec:factorization}
In this section, we analyze how maps between \textit{F}-finite animated $\FF_p$-algebras can be decomposed into simpler maps. As a first step, we consider the case of relatively semiperfect maps.

\begin{con}
\label{con:factorization_f-surjective_map_der}
Let $R\to S$ be a relatively semiperfect map of \textit{F}-finite animated $\FF_p$-algebras. We choose a free $\FF_p$-algebra of finite type along with a relatively semiperfect map
\[
A=\FF_p[x_1,\dots,x_n] \to R
\]
for some $n\ge0$. Then $R\to S$ becomes a map in $\sperf{A}$.
Since $F_*A$ is a finite free $A$-module, the relative inverse limit perfection over $A$ is relatively perfect and thus
$R^{\perfection/A}\to S^{\perfection/A}$ is a map in $\perf{A}$. The counit $\varepsilon$ of the adjunction from Proposition \ref{prop:right_adjoint_animated} yields the following commutative diagram in the category $\sperf{A}$:
\[
\begin{tikzcd}
{R^{\perfection/A}} \arrow[r] \arrow[d, "\varepsilon_R"] & {S^{\perfection/A}} \arrow[d, "\varepsilon_S"] \\
{R} \arrow[r]           & {S}          
\end{tikzcd}
\]
This yields a factorization of the map $R\to S$:
\[
R \to R\otimes^L_{R^{\perfection/A}}S^{\perfection/A} \to S.
\]
The first of these maps is relatively perfect, and the second is surjective on $H^0$. Indeed, since $R^{\perfection/A}\to S^{\perfection/A}$ is a map in $\perf{A}$, both the map itself and its base-change are relatively perfect by Lemma \ref{lem:relative_perfect_basic_properties_der}. Moreover, the second map is surjective on $H^0$, as the same holds for $S^{\perfection/A}\to S$ by Lemma~\ref{lem:connective_surjective_der}.
\end{con}


\begin{lem}
\label{lem:f-srujective_factorization_der}
Let $R\to S$ be a relatively semiperfect map of \textit{F}-finite  animated $\FF_p$-algebras. Construction \ref{con:factorization_f-surjective_map_der} gives a factorization
\[
R \to T \to S
\]
such that $T$ is relatively perfect over $R$, and the second map is surjective on $H^0$. Moreover, we have:
\begin{enumerate}
\item[(i)] If $H^0(R)$ is Noetherian and $R$ is $n$-coconnective, then $T$ is also $n$-coconnective.
\item[(ii)] If both $R$ and $H^0(S)$ are Noetherian, then $T$ is Noetherian as well.
\end{enumerate}
\end{lem}
\begin{proof}
We use the same notation as in Construction \ref{con:factorization_f-surjective_map_der}. If $H^0(R)$ is Noetherian, then $R^{\perfection/A}$ is also Noetherian by Proposition \ref{prop:perfection_is_Noetherian_der} and therefore regular. Hence, by \citep[Proposition 5.2]{Kat86}, $S^{\perfection/A}$ is flat over $R^{\perfection/A}$, which implies that $T=R\otimes^L_{R^{\perfection/A}}S^{\perfection/A}$ is flat over $R$. 
Consequently, 
we obtain canonical isomorphisms for all $k\le 0$:
\[
H^k(T)
\cong 
H^k(R)\otimes_{R^{\perfection/A}} S^{\perfection/A} 
\cong 
H^k(R)\otimes_{H^0(R)} H^0(T).
\]
The first isomorphism establishes statement (i). Now, assume that $R$ and $H^0(S)$ are Noetherian. By the same arguments as above, $S^{\perfection/A}$ is Noetherian. Therefore, by the first isomorphism, $H^0(T)$ is a quotient of a Noetherian ring. The second isomorphism shows, that the negative cohomology groups of $T$ are finitely generated over $H^0(T)$, as the same holds for $R$ by assumption.
\end{proof}

By Proposition \ref{prop:cotangent_rel_perf_animated}, the cotangent complex of any relatively perfect map is acyclic. For maps of Noetherian \textit{F}-finite rings, the converse also holds. 


\begin{thm}
\label{thm:acyclic_cotanget_implies_rel_perf}
Let $ R \to S $ be a map of Noetherian $ F $-finite rings. If $ L_{S/R} \simeq 0 $, then $ S $ is relatively perfect over $ R $.
\end{thm}
\begin{proof}
By assumption, $\Omega_{S/R}=0$, so $S$ is relatively semiperfect over $R$ by Lemma  \ref{lem:cotangent_semiperfect_der}. The previous lemma now yields a factorization
\[
R\to T \to S,
\]
where $T$ is Noetherian, the map $R\to T$ is relatively perfect, and the map $T\to S$ is surjective. It suffices to show that $S$ is relatively perfect over $T$. 

The canonical triangle of cotangent complexes shows that $L_{S/T}\simeq 0$, as both $L_{S/R}$ and $L_{T/R}$ vanish by assumption and by Proposition~\ref{prop:cotangent_rel_perf_animated}, respectively. In particular, $T\to S$ is \'etale and thus relatively perfect by \cite[\href{https://stacks.math.columbia.edu/tag/0EBS}{Lemma 0EBS}]{stacks-project}. 
\end{proof}



\begin{cor}
Let $R$ be a discrete Noetherian \textit{F}-finite ring and $I \subset R$ an ideal. Then the $I$-adic completion $R^{\wedge}_I$ is a relatively perfect $R$-algebra. In other words, the discrete relative inverse limit perfection of $R/I$ over $R$ is relatively perfect.
\end{cor}
\begin{proof}
By \citep[Proposition 3.11]{BF25}, the cotangent complex of the natural map $R\to R^{\wedge}_I$ vanishes. Thus, by Theorem~\ref{thm:acyclic_cotanget_implies_rel_perf}, the map is relatively perfect.
\end{proof}

The following corollary shows, for instance, that a Noetherian $F$-finite $\FF_p$-algebra is regular with a $p$-basis in the sense of \citep[Définition 21.1.9]{EGAIV} if and only if it is relatively perfect over polynomial ring over $\FF_p$ in finitely many variables.

\begin{cor}
\label{cor:p-basis}
Let $R\to S$ be a map of discrete Noetherian \textit{F}-finite rings. Then $S$ admits a $p$-basis of length $n\ge0$ over $R$, and $L_{S/R}\simeq \Omega^1_{S/R}[0]$, if and only if $S$ is relatively perfect over a polynomial $R$-algebra in $n$ variables.
\end{cor}
\begin{proof}
The "if" direction follows from Proposition \ref{prop:cotangent_rel_perf_animated} since the cotangent complex of a polynomial algebra is concentrated in degree $0$, and the variables form a $p$-basis. Conversely, 
let ${s_1,\dots,s_n}$ be a $p$-basis of $S$ over $R$ and define a map $R[x_1,\dots,x_n] \to S$ by sending each free variable $x_i$ to $s_i$. Since, by \citep[Corollaire 21.2.5]{EGAIV}, any $p$-basis is also a differential basis, the induced map $L_{R[x_1,\dots,x_n]/R}\otimes^L_{R[x_1,\dots,x_n]}S \to L_{S/R}$ is an equivalence. Therefore, the canonical fiber sequence of cotangent complexes associated with $R\to R[x_1,\dots,x_n] \to S$ shows that $L_{S/R[x_1,\dots,x_n]}$ vanishes. By Theorem \ref{thm:acyclic_cotanget_implies_rel_perf}, this implies that $S$ is relatively perfect over $R[x_1,\dots,x_n]$.
\end{proof}

For a formally smooth map of locally Noetherian $F$-finite $\FF_p$-schemes we obtain a factorization result similar to \cite[\href{https://stacks.math.columbia.edu/tag/039Q}{Theorem 039Q}]{stacks-project}.

\begin{cor}
\label{cor:factorizationofschemes}
Let $f:X\to Y$ be a formally smooth morphism of locally Noetherian $F$-finite schemes. Then for any $x\in X$ and any affine open neighborhood $V\subset Y$ of $f(x)$, there exists an integer $n\ge0$, an affine open neighborhood $U$ of $x$ with $f(U)\subset V$, and a commutative diagram
\[
\begin{tikzcd}
X \arrow[d] & U \arrow[l] \arrow[d] \arrow[r, "g"] & \mathbf{A}^n_V \arrow[d]                \\
Y           & V \arrow[l]                     & V \arrow[l, equal]
\end{tikzcd}
\]
in which $g$ is formally \'etale.
\end{cor}
\begin{proof}
Let $x\in X$, and let $V\subset Y$ an affine open neighborhood of $f(x)$. By choosing an affine open subscheme $W\subset X$ containing $x$ with $f(W)\subset V$, we may reduce to the case where $X$ and $Y$ are affine; see \cite[\href{https://stacks.math.columbia.edu/tag/02H3}{Lemma 02H3}]{stacks-project}. In this case, $f$ corresponds to a formally smooth map $R \to S$ of Noetherian \textit{F}-finite $\mathbb{F}_p$-algebras. By \citep[Proposition 2.15]{BF25}, the cotangent complex $L_{S/R}$ is concentrated in degree zero and $\Omega_{S/R}$ is a finite projective $S$-module. In particular, $\Omega_{S/R}\otimes_{S}\kappa(x)$ is an $n$-dimensional $\kappa(x)$-vector space for some integer $n\ge 0$. Let $d(s_1),\dots, d(s_n)\in \Omega_{S/R}$ be elements whose images under the canonical map form a basis of this vector space. These elements induce a map $S^n\to \Omega_{S/R}$ which becomes an isomorphism on some standard open neighborhood $U=D(a)$ of $x$:
\[
S^n_a \to[\cong] (\Omega_{S/R})_a \cong \Omega_{S_a/R}.
\]
Thus, $S_a$ admits a differential basis over $R$. Since, by \citep[Theorem 1]{Tyc88}, a differential basis of Noetherian \textit{F}-finite rings is also a $p$-basis, the claim follows from the previous corollary.
\end{proof}

Applying Lemma \ref{lem:f-srujective_factorization_der}, we obtain the following result, showing that any map of \textit{F}-finite animated rings admits a free-by-surjective factorization up to a relatively perfect map.

\begin{thm}
\label{thm:F-finite_factorization_der}
Let $R \to S$ be a map of \textit{F}-finite animated $\FF_p$-algebras, and let $R\to R^{\prime} \to S$ be a factorization through a free $R$-algebra of finite type such that the second map is relatively semiperfect. Applying \ref{con:factorization_f-surjective_map_der} to the relatively semiperfect map yields a factorization
\[
R \to R^{\prime} \to T \to S
\]
where the first map is free of finite type, the second map is relatively perfect and the third map is surjective on $H^0$. Moreover:
\begin{enumerate}
\item[(i)] If $H^0(R)$ is Noetherian and $R$ is $n$-coconnective, then both $R'$ and $T$ are $n$-coconnective.
\item[(ii)] If both $R$ and $H^0(S)$ are Noetherian, then all animated rings appearing in the factorization are Noetherian.
\end{enumerate}
\end{thm}
\begin{proof}
Since $R^{\prime}$  is flat over $R$ it is $n$-coconnective whenever $R$ is $n$-coconnective. By Hilbert's basis theorem, $H^0(R^{\prime})$ is Noetherian if $H^0(R)$ has this property. Moreover, if $R$ is Noetherian, the flatness of $R^{\prime}$ over $R$ ensures that $R^{\prime}$ is Noetherian as well. Thus, the conclusion follows directly from Lemma \ref{lem:f-srujective_factorization_der} applied to $R^{\prime}\to S$.
\end{proof}

In the case where the rings are Noetherian and discrete, the following provides a new proof of \citep[Theorem 4.4]{BF25}.

\begin{cor}
\label{cor:factorization_discrete_der}
Let $R\to S$ be a homomorphism of ordinary \textit{F}-finite $\FF_p$-algebras. Theorem \ref{thm:F-finite_factorization_der} gives a factorization
\[
R \to R' \to T \to S
\]
where the first map is free of finite type, the second map is relatively perfect and the third map is surjective. If $R$ is Noetherian, then $T$ is discrete, and if both $R$ and $S$ are Noetherian, then $T$ is Noetherian as well.
\end{cor}
\begin{proof}
This follows from the previous theorem: If $R$ is Noetherian and discrete, then every ring in the factorization is discrete. Moreover, if both $R$ and $S$ are Noetherian, then all rings in the factorization are Noetherian as well.
\end{proof}

\begin{rem}
\label{rem:final_remark}
If $R$ in the previous corollary is non-Noetherian, the animated ring $T$ may not be discrete; see Construction \ref{con:factorization_f-surjective_map_der}. Nevertheless, we still obtain a factorization of classical rings
\[
R \to R^{\prime} \to H^0(T) \to S,
\]
and the discrete relative Frobenius of $H^0(T)$ over $R^{\prime}$ is an isomorphism. In particular, $H^0(T)$ is formally \'etale over $R^{\prime}$ by \citep[Théorème 21.2.7]{EGAIV}.

However, by the first part of Lemma~\ref{lem:relative_perfect_basic_properties_der} and  Corollary \ref{cor:equivalence_animated}, the map $R^{\prime} \to H^0(T)$ is not relatively perfect  in the sense of our definition unless $T$ itself is discrete; see also Remark \ref{rem:tor_independence_relative_perfect_der}.
\end{rem}


\emergencystretch=1em
\printbibliography
\end{document}